\documentclass[11pt]{amsart}
\usepackage[sorted-cites]{amsrefs}
\usepackage{verbatim}
\usepackage{breqn}
\usepackage{color}
\usepackage{amscd}
\usepackage{amsmath,amsthm,amscd}
\usepackage{enumitem}
\usepackage{pgf,tikz}

\newfont{\fnt}{cmr10 scaled 550}

\newtheorem{theorem}{Theorem}
\newtheorem*{theorem*}{Theorem}

\newtheorem{lemma}{Lemma}
\newtheorem{cor}{Corollary}

\theoremstyle{remark}
\newtheorem{remark}{Remark}
\newenvironment{theproof}{\noindent{\textbf{Proof: }}}{ \begin{flushright}$\square$\end{flushright}}
    \def\proof{\vspace{12pt}\begin{theproof}}
    \def\proofend{\end{theproof}}

\renewcommand{\epsilon}{\varepsilon}

\newcommand{\ric}{\mathop{\mathrm{Ric}}\nolimits}
\newcommand{\hess}{\mathop{\mathrm{Hess}}\nolimits}
\newcommand{\vol}{\mathop{\mathrm{vol}}\nolimits}
\newcommand{\dom}{\mathop{\mathrm{dom}}\nolimits}
\newcommand{\supp}{\mathop{\mathrm{supp}}\nolimits}
\newcommand{\inner}[1]{\left\langle#1\right\rangle}
\newcommand{\innerS}[1]{\langle#1\rangle}
\newcommand{\innerP}[1]{\left(#1\right)}
\newcommand{\norm}[1]{\left\|#1\right\|}
\newcommand{\module}[1]{\left|#1\right|}
\newcommand{\lnorm}[3]{\norm{#1}_{L^{#2}(#3)}}
\newcommand{\lnormm}[2]{\norm{#1}_{L^{#2}(M)}}
\newcommand{\lfnorm}[3]{\norm{#1}_{L^{#2}_f(#3)}}
\newcommand{\lfnormm}[2]{\norm{#1}_{L^{#2}_f(M)}}
\newcommand{\branch}[1]{\left\{\begin{array}{ll}#1\end{array}\right.}
\newcommand{\dd}[2]{\frac{\partial #1}{\partial #2}}
\newcommand{\rtilde}{\tilde{r}}
\newcommand{\R}{\mathbb{R}}
\newcommand{\N}{\mathbb{N}}
\def\grad{\nabla}

\begin{document}

    \title[On the essential spectrum of the Laplacians]
    {On the essential spectrum of the Laplacian and the drifted Laplacian}
    \date{January 29, 2013}
    \author[Leonardo Silvares]{Leonardo Silvares}
    \address{Departamento de Matem\'atica e Estat\'istica, UNIRIO - Universidade Federal do Estado do Rio de Janeiro, Av. Pasteur 458, Rio de Janeiro, Brazil}
    \email{leo.silvares@uniriotec.br}

    \begin{abstract}
        This paper concerns the $L^2$ essential spectrum of the Laplacian $\Delta$ and the drifted Laplacian $\Delta_f$ on complete Riemannian manifolds endowed with a weighted measure $e^{-f}d\vol_g$. We prove that the essential spectrum of the drifted Laplacian $\Delta_f$ is $[0,\infty) $ provided the Bakry-\'Emery curvature tensor $\ric_f$ is nonnegative and $f$ has sublinear growth. When  $\ric_f \geq \frac{1}{2}g$ and $|\nabla f|^2 \leq f$, we show that the essential spectrum of the Laplacian is also $[0,\infty)$. During the proofs of these results, the $f$-volume growth estimate plays an important role and may be of independent interest.
    \end{abstract}

\maketitle

\section{Introduction}

    The spectra of the Laplacian of complete non-compact manifolds has been studied and computed for a large class of manifolds in the past two decades. When the manifold has a soul whose exponential map is a diffeomorphism, supposing non-negative sectional curvature and some other additional conditions, Escobar \cite{escobar} and Escobar-Freire \cite{escobar-freire} proved that the $L^2$ spectrum of the Laplacian is $[0, \infty)$. The additional conditions assumed by them were proved to be unnecessary by Zhou \cite{zhou}.

    If the manifold has a pole, the essential spectrum of the Laplacian was proved to be $[0, \infty)$ by Li \cite{li} (supposing the Ricci curvature non-negative) and by Chen-Lu \cite{chen-lu} (supposing the radial sectional curvature non-negative). The same result was proved by Donnelly \cite{donnelly} when the Ricci curvature is non-negative and the volume has Euclidean growth.

    Without having to assume the existence of a pole, J. Wang \cite{wang1} proved that the essential spectrum $L^p$ of the Laplacian is $[0, \infty)$, for $p \in [1, \infty)$, provided the Ricci curvature is bounded below by $-\delta/r^2$, where $r$ is the distance from a fixed point on the manifold and $\delta > 0$ depends only on the dimension. This result was generalized by Lu and Zhou \cite{lu-zhou}, assuming only that $\underline{\lim}_{r\rightarrow\infty}\ric = 0$. Both works proved initially that the $L^1$  essential spectrum is $[0, \infty)$ and extended this fact to $L^p$, $p \in [1, \infty]$, by using a result proved by Sturm \cite{sturm} which asserted the $p$-independence of the spectrum when the volume has uniform sub-exponential growth and the Ricci curvature is bounded below.

    Circumventing the use of Sturm's Theorem \cite{sturm} by a generalization of Classical Weyl's Criterion, Charalambous and Lu \cite{charalambous1} extended the results previously obtained in $L^2$ to the case where $\ric$ is asymptotically non-negative and $\overline{\lim}_{r\rightarrow\infty} \Delta r \leq 0$. 

    In many occasions, it is natural to consider a manifold $M$ endowed with a weighted measure of the form $e^{-f}d\vol_g$, where $f$ is a smooth function. We have the drifted Laplacian $\Delta_f$, defined by
    \[\Delta_f u = \Delta u - \inner{\grad f, \grad u},\]
    associated with the differential of the volume form of $e^{-f}d\vol_g$ the same way $\Delta$ is associated to $d\vol_g$. Moreover, $\Delta_f$ is a self-adjoint operator on the space $L^2_f$ of square integrable functions on $M$ with respect to the measure $e^{-f}d\vol_g$. 
    A natural extension of the Ricci curvature tensor to this new context is the Bakry-\'Emery Ricci curvature tensor, given by
    \[\ric_f = \ric + \hess f.\]

    Manifolds with the Bakry-\'Emery Ricci curvature bounded below have been studied by many authors in recent years, with the special objective of obtaining similar results to the case where the Ricci curvature is bounded below. Many results have been obtained, particulary some interesting \linebreak volume estimates and splitting theorems, which can be found in \cite{wei1} and \cite{munteanu}. In this paper, we apply some of these results to the analysis of the drifted Laplacian $\Delta_f$. A natural context for the Bakry-\'Emery Ricci curvature tensor is the study of \textit{gradient Ricci solitons}, manifolds where $\ric_f = \lambda g$, for some constant $\lambda$ and function $f$.

    Our main goal in this paper is to determine the $L^2$ and $L^2_f$ essential spectrums of $-\Delta$ and $-\Delta_f$ respectively, supposing some lower bounds for $\ric_f$ and imposing conditions on $f$.

    We state now our first result. Fixing $p \in M$ and denoting $r(x) = d(p, x)$, we have:

    \begin{theorem}\label{principal1ext}
         Let $M$ be a non-compact complete manifold. If $f: M \rightarrow \R$ is a smooth function such that $\ric_f \geq 0$ and $\lim_{r \rightarrow \infty}\frac{|f|}{r} = 0$, the $L^2_f(M)$ essential spectrum of $-\Delta_f$ is $[0, \infty).$
    \end{theorem}

    \begin{remark}
    In order to prove this Theorem, we could suppose only \linebreak $\ric_f(\grad r, \grad r)$, provided we assume, in the case $\vol_f(M) < \infty$, that \linebreak $\vol_f(B_p(r))$ does not decay exponentially at $p$. We are denoting by $\vol_f$ the volume calculated using the density $e^{-f}d\vol_g$.
    \end{remark}

    \begin{remark}
    We remark that the hypotheses $\lim_{r \rightarrow \infty}\frac{|f|}{r} = 0$ of Theorem \ref{principal1ext} is somehow sharp. For $a > 0$, it is known an example of $M$ and $f$ such that $\ric_f \geq 0$, $\lim_{r \rightarrow \infty}\frac{|f|}{r} = a$ and the $L^2_f(M)$ essential spectrum of $-\Delta_f$ is not $[0, \infty)$.
    Consider a gradient steady Ricci soliton, i.e. a manifold $(M, g)$ satisfying $\ric_f = R_{ij} + f_{ij} = 0$, and take $a > 0$ such that $|\grad f|^2 + S = a^2$, where $S$ is the scalar curvature of $M$. It is well known that $f$ has linear growth and, by \cite{munteanu2}, the $L^2$ spectrum of $-\Delta_f$ is contained in $[\frac{1}{4}a^2, \infty)$.
    \end{remark}

    The following theorems, which are used in the proof of Theorem \ref{principal1ext}, are also of independent interest.

    \begin{theorem}\label{controle_crescimento_teorema1ext}
        Let $M$ be a non-compact complete manifold. If $f: M \rightarrow \R$ is a smooth function such that $\ric_f \geq 0$ and $\lim_{r \rightarrow \infty}\frac{|f|}{r} = 0$,
        then \[\lim_{r \rightarrow \infty}\frac{\vol_f(B_p(r))}{e^{Cr}} = 0,\]
        for all $C > 0$.
    \end{theorem}
    \begin{theorem}\label{decaimento_nao_exponencial_teorema1ext}
        Let $M$ be a non-compact complete manifold. If $f: M \rightarrow \R$ is a smooth function such that $\ric_f \geq 0$, $\lim_{r \rightarrow \infty}\frac{|f|}{r} = 0$ and $\vol_f(M) < \infty$, then the $f$-volume of $M$ does not decay exponentially at $p$.
    \end{theorem}

    One may naturally ask what happens to essential spectrum of $-\Delta_f$ when the condition $\ric_f \geq 0$ is strengthened to $\ric_f \geq \frac{1}{2}g$. 
    According to \cite{naber}, it is empty since the spectrum of $-\Delta_f$ is discrete. However, we can make te same question about $-\Delta$. In the case $\ric_f = \frac{1}{2}g$, Lu and Zhou \cite{lu-zhou} proved that the essential spectrum of $-\Delta$ is $[0,\infty)$, supposing uniform sub-exponencial volume growth or $\lim R/r^2 = 0$. Later, Charalambous and Lu \cite{charalambous1} proved this result without such hypotheses. Notice that this is the case of a  gradient shrinking Ricci soliton, because its equation
    \[R_{ij} + \grad_i\grad_j f = \rho g_{ij}, \rho > 0,\]
    can be normalized to $\rho = \frac{1}{2}$. 

    We have the following theorem:

    \begin{theorem}\label{principal2}
         Let $M$ be a non-compact complete manifold. If $f: M \rightarrow \R$ is a smooth function such that $\ric_f \geq \frac{1}{2} g$ and $|\grad f|^2\leq f$, the $L^2(M)$ essential spectrum of $-\Delta$ is $[0, \infty).$
    \end{theorem}

    This paper is organized as follows. In Section \ref{sec_weyl}, we start with a brief review of the Weyl's Criterion for self-adjoint operators in a Hilbert space, and adapt this criterion for the particular use in $L^2$ and $L_f^2$. In Section \nolinebreak \ref{rtil_teorema1}, we present a technique for obtaining smooth approximations of distance functions, and use the hypotheses on $\ric_f$ and $f$ in order to ensure some useful estimates on the decay of the approximate distances. In Section \ref{estimativas_vol_f_teorema1ext} we present a sequence of Lemmas and volume estimates Theorems which will lead to the proof of the first main Theorem, which will be done in Section \ref{prova_do_teorema1ext}. In Section \ref{prova_do_teorema2}, assuming the hypotheses in Theorem \ref{principal2}, we adapt the estimates obtained in Section \ref{estimativas_vol_f_teorema1ext} and prove Theorem \ref{principal2}.\\

    \noindent \textbf{Acknowledgement:} The author would like to thank his PhD advisor at Universidade Federal Fluminense (Brazil), Professor Detang Zhou, whose support, help and teachings have been essential to the author's mathematical and personal development. The author also wants to thank Guofang Wei for her interest on this work and some remarks.

\section{Notations and Weyl's Criteria}\label{sec_weyl}

    Let $H$ be a self-adjoint operator in a Hilbert space $\mathcal{H}$ whose inner product is denoted by $(\cdot,\cdot)$ and its norm by $\norm\cdot$.
    We denote by $\sigma(H)$ the spectrum of $H$ and by $\sigma_{\mathrm{ess}}(H)$ its essential spectrum.

    \begin{theorem}
        (Classical Weyl's Criterion) A nonnegative real number $\lambda$ belongs to $\sigma(H)$ if, and only if, there exists a sequence $\{u_i\}_{i \in \N} \subset \dom(H)$ such that

        \begin{enumerate}[label={(\arabic*)}]
            \item $\norm{u_i} = 1, \forall i$;
            \item $(H - \lambda)u_i \rightarrow 0$ in $\mathcal{H}$, as $i \rightarrow \infty$.
        \end{enumerate}
        Moreover, $\lambda$ belongs to $\sigma_{\mathrm{ess}}(H)$ if, and only if, there exists a sequence \linebreak $\{u_i\}_{i \in \N} \subset \dom(H)$ satisfying (1), (2) and
        \begin{enumerate}[label={(\arabic*)}]
            \setcounter{enumi}{2}
            \item $u_i \rightarrow 0$ weakly in $\mathcal{H}$.\\
        \end{enumerate}
    \end{theorem}

    As a consequence of the Weyl's Criterion, we have the following Theorem, which is consequence of a more general version in \cite{charalambous1}:

    \begin{theorem}\label{weyl_generalizado}
        A nonnegative real number $\lambda$ belongs to $\sigma(H)$ if, and only if, there exists a sequence $\{u_i\}_{i \in \N} \subset \dom(H)$ such that
        \begin{enumerate}[label={(\arabic*)}]
            \item $\norm{u_i} = 1, \forall i$;
            \item $\left((H + 1)^{-1}u_i, (H - \lambda)u_i\right) \rightarrow 0$ as $i \rightarrow \infty$.
            \item $\left(u_i, (H - \lambda)u_i\right) \rightarrow 0$ as $i \rightarrow \infty$.
        \end{enumerate}
        Moreover, $\lambda$ belongs to $\sigma_{\mathrm{ess}}(H)$ if, and only if, there exists a sequence \linebreak $\{u_i\}_{i \in \N} \subset D(H)$ satisfying (1) - (3) and
        \begin{enumerate}[label={(\arabic*)}]
            \setcounter{enumi}{3}
            \item $u_i \rightarrow 0$ weakly in $\mathcal{H}$.\\
        \end{enumerate}
    \end{theorem}

    Given $R \subset M$, $p \in [1, \infty)$ and $f \in C^\infty(M)$, we define \linebreak $L^p(R):= L^p(R, g, d\vol_g)$, $L^p_f(R) := L^p(R, g, e^{-f}d\vol_g)$ and the norms
    \begin{gather*}
        \norm{u}_{L^p(R)} := \norm{u}_{L^p(R,g,d\vol_g)}  := \left(\int_R |u|^p\;d\vol_g \right)^{\frac{1}{p}},\\
        \norm{u}_{L_f^p(R)} := \norm{u}_{L^p(R,g,e^{-f}d\vol_g)} := \left(\int_R |u|^p\;e^{-f}d\vol_g \right)^{\frac{1}{p}}.
    \end{gather*}

    Sometimes we write simply $L^p = L^p(M)$ and $L^p_f = L^p_f(M)$.

    Notice that, endowed with the inner products $(u,v) = \int_M u\overline{v} d\vol_g$ and $(u,v)_f = \int_M u\overline{v} e^{-f}d\vol_g$, respectively, $L^2$ and $L^2_f$ are Hilbert spaces in which $-\Delta$ and $-\Delta_f$ are self-adjoint (unbounded) operators.

    Defining
    \[\lnorm{u}{\infty}{R} := \sup_{x \in R}\{|u(x)|\}, \lnormm{u}{\infty} := \sup_{x \in M}\{|u(x)|\},\]
    and denoting by $C^\infty_0(M)$ the set of compactly supported functions $M \rightarrow \R$, we have the following result:\\

    \begin{cor}\label{weyl_lap_f}
        If there exists a sequence $\{u_i\} \subset C^\infty_0(M)$ such that
        \begin{enumerate}
            \item $\displaystyle \frac{\lnormm{u_i}{\infty}\;\lfnorm{(-\Delta_f - \lambda)u_i}{1}{M}}{\lfnormm{u_i}{2}^2} \rightarrow 0$ as $i \rightarrow \infty$;
            \item for any compact $K \subset M$, there exists $i_0$ such that the support of $u_i$ is outside $K$ for $i > i_0$; and
            \item $\partial(\supp(u_i))$ is a $C^\infty$ $(n-1)$-submanifold of $M$,\\
        \end{enumerate}
        then $\lambda \in \sigma_{\mathrm{ess}}(-\Delta_f)$ .
    \end{cor}

    \proof
        Without loss of generality, we assume that $\lfnormm{u_i}{2} = 1$ for all $i$. Hence, \[\lnormm{u_i}{\infty}\;\lfnormm{(-\Delta_f - \lambda)u_i}{1} \rightarrow 0.\]

        Since
        \begin{eqnarray*}
            \module{\innerP{u_i, (-\Delta_f-\lambda)\;u_i}_f}&=&     \int_M u_i\,\cdot\,(-\Delta_f-\lambda)u_i\;e^{-f}d\vol_g\\
                                                            &\leq&  \int_M \module{u_i}\,\cdot\,\module{(-\Delta_f-\lambda)u_i}\;e^{-f}d\vol_g\\
                                                            &\leq&  \lnormm{u_i}{\infty}\,\cdot\,\int_M \module{(-\Delta_f-\lambda)u_i}\;e^{-f}d\vol_g\\
                                                            &=&     \lnormm{u_i}{\infty}\,\cdot\,\lfnormm{(-\Delta_f-\lambda)u_i}{1},
        \end{eqnarray*}
        the condition (3) of Theorem \ref{weyl_generalizado} is satisfied. We shall now verify \linebreak condition (2).

        Since each $D_i := \supp(u_i)$ is a compact set with smooth boundary and $\Delta_f$ is an elliptic self-adjoint operator, there exists $v_i \in C^\infty(D_i)$ such that
        \begin{gather*}
            \left.((-\Delta_f + 1)v_i)\right|_{D_i} = u_i,\\
            v_i|_{\partial D_i} = 0.
        \end{gather*}

        For each $i$, $D_i$ is a compact set, so we take $x_i \in D_i$ such that \linebreak $|v_i(x_i)| = \sup_{x \in D_i}|v_i(x)| = \lnorm{v_i}{\infty}{D_n}$. Without loss of generality, we suppose $v_i(x_i) > 0$. Therefore
        \[\Delta v_i(x_i) \leq 0, \grad v_i(x_i) = 0.\]

        Now we have
        \begin{eqnarray*}
            \|u_i\|_{L^\infty}
                &=&     \|u_i\|_{L^\infty(D_i)}\\
                &=&     \|(-\Delta_f + 1)v_i\|_{L^\infty(D_i)}\\
                &\geq&  |(-\Delta_f + 1)v_i)(x_i)| \\
                &=&     |- (\Delta_f v_i)(x_i) + v_i(x_i)|\\
                &=&     |- (\Delta v_i)(x_i) + \inner{\grad f, \grad v_i}(x_i) + v_i(x_i)|\\
                &\geq& v_i(x_i)\\
                &=& \|v_i\|_{L^\infty},
        \end{eqnarray*}
        which leads to

        $\displaystyle\module{\innerP{(-\Delta_f+1)^{-1}u_i, (-\Delta_f-\lambda)u_i}_f} = $
        \begin{eqnarray*}
            \phantom{\module{\innerP{(-\Delta_f+1)^{-1}u_i, }}}
                    &=&     \int_M ((-\Delta_f+1)^{-1}u_i)\,\cdot\,(-\Delta_f-\lambda)u_i\;e^{-f}d\vol_g\\
                    &=&     \int_{D_i} ((-\Delta_f+1)^{-1}u_i)\,\cdot\,(-\Delta_f-\lambda)u_i\;e^{-f}d\vol_g\\
                    &\leq&  \int_{D_i} \module{(-\Delta_f+1)^{-1}u_i}\,\cdot\,\module{(-\Delta_f-\lambda)u_i}\;e^{-f}d\vol_g\\
                    &\leq&  \int_{D_i} \module{v_i}\,\cdot\,\module{(-\Delta_f-\lambda)u_i}\;e^{-f}d\vol_g\\
                    &\leq&  \lnorm{v_i}{\infty}{D_i}\,\cdot\,\int_M \module{(-\Delta_f-\lambda)u_i}\;e^{-f}d\vol_g\\
                    &\leq&     \lnormm{u_i}{\infty}\,\cdot\,\lfnormm{(-\Delta_f-\lambda)u_i}{1}.
        \end{eqnarray*}

        This proves that $\module{\innerP{(-\Delta_f+1)^{-1}u_i, (-\Delta_f-\lambda)u_i}_f} \rightarrow 0$.
    \proofend

    Analogous to previous Corollary, we have

    \begin{cor}\label{weyl_lap}
        If there exists a sequence $\{u_i\} \subset C^\infty_0(M)$ such that
        \begin{enumerate}
            \item $\displaystyle \frac{\lnormm{u_i}{\infty}\;\lnormm{(-\Delta_f - \lambda)u_i}{1}}{\lnormm{u_i}{2}^2} \rightarrow 0$ as $i \rightarrow \infty$;
            \item For any compact $K \subset M$, there exists $i_0$ such that the support of $u_i$ is outside $K$ for $i > i_0$; and
            \item $\partial(\supp (u_i))$ is a $C^\infty$ $(n-1)$-submanifold of $M$,\\
        \end{enumerate}
        then $\lambda \in \sigma_{\mathrm{ess}}(-\Delta)$ .
    \end{cor}

\section{$C^\infty$ approximations of $r$}\label{rtil_teorema1}

     In this section, we fix $p \in M$ and $f \in C^\infty(M)$ and let $r(x) = d(x,p)$. Given a decreasing continuous function $\delta: \R^+ \rightarrow \R^+$, such that \linebreak $\lim_{r \rightarrow \infty} \delta(r) = 0$, we are going to construct functions $b$ and $\rtilde$ in $C^\infty(M)$ satisfying

    \begin{enumerate}[label={\roman*.}]
        \item $\displaystyle \norm{b}_{L^1_f(M \setminus B_p(r))} \leq \delta(r)$;
        \item $\displaystyle \norm{\grad\rtilde - \grad r}_{L^1_f(M \setminus B_r(p))} \leq \delta(r)$;
        \item $\displaystyle \module{\rtilde(x) - r(x)} \leq \delta(r(x))$ and $\module{\grad\rtilde(x)}\leq 2, \forall x \in M, r(x) > 2$, and
        \item $\displaystyle \Delta_f\rtilde(x) \leq \sup_{y \in B_x(1)}\{\Delta_f r(y)\} + \delta(r(x)) + \module{b(x)}, \forall x \in M, r(x) > 2$ in the sense of distributions.
    \end{enumerate}

    \begin{remark}
    We would like to emphasize the fact that the functions $\rtilde$ and $b$ depend on the given $\delta$. In the proof of forthcoming lemmas and theorems, some conditions on the decay of $\delta$ will be imposed; we will consider $\delta$ satisfying all those conditions.
    \end{remark}

    The construction of such an approximation is done in \cite{lu-zhou} and \cite{charalambous1} for $L^1$ and $\Delta$ instead of $L^1_f$ and $\Delta_f$. The adaptation is not difficult, but, for the sake of completeness, we present the construction here.

    We shall now obtain a suitable partition of the unity. Taking $x \in M$, we consider a neighborhood $U_x$ such that

    \begin{itemize}
        \item there exist local coordinates $\mathbf{x}_i$ defined in an open set containing $\overline{U_x}$ such that $x = (0, \ldots, 0)$; and\\
        \item if $\grad f (x) = f_i(x) \partial x^i$, then $|\grad f(y) - \sum_i f_i(x)\partial x^i (y)| < \delta(r(y))/4$ for all $y \in U_x$.
    \end{itemize}

    Moreover, we reduce the $U_x$ in such way that $U_x \subset B_x(1)$ for all $x$.

    We may extract a locally finite cover $\{U_i\}_{i \in \N}$ from $\{U_x\}_{x \in M}$, and let $\{\psi_i\}$ be the associated partition of the unity.

    Let $\xi(\mathbf{x})$ be a non-negative smooth function whose support is within the unit ball of $\R^n$ and such that
    \[\int_{\R^n}\xi d\mathbf{x} = 1.\]

    Fixing $i \in \N$ and $\epsilon > 0$, we denote $\mathbf{x}_i = (x_i^1, ... x_i^n)$ the local coordinates of $U_i$ (as above). Define
    \begin{equation}\label{r_i}
        r_{i,\epsilon_i}(\textbf{x}_i) = \frac{1}{\epsilon_i^n}\int_{\R^n} \xi\left(\frac{\textbf{y}_i}{\epsilon_i}\right)r(\textbf{x}_i + \textbf{y}_i)\;d\textbf{y}_i
    \end{equation}
    and $m_i := \sup_{U_i}\{ |\grad \psi_i| + |\Delta_f \psi_i| \}$. We take $\epsilon_i$ em (\ref{r_i}) such that, for all $x \in U_i$,
    \begin{gather}
        \label{controle_r_menos_r_i}
        |r_{i,\epsilon_i}(x) - r(x)| \leq \frac{\delta(r(x))}{2^im_i},\\
        \label{controle_grad_r_i_pontual}
        |\grad r_{i,\epsilon_i}(x)| < 2
    \end{gather}
    and
    \begin{equation}\label{controle_grad_r_menos_grad_r_i}
        \left\|\grad r_{i,\epsilon_i}(x) - \grad r(x)\right\|_{L^1_f(U_i)} \leq \frac{\delta(r(x)-1)}{2^im_i}.
    \end{equation}

    Define, for $x \in M$,
    \[\rtilde(x) = \sum_i \psi_i(x) r_{i,\epsilon_i}(x)\;\;\textnormal{ and }\;\; b(x) = 2\sum_i\inner{\grad \psi(x), \grad r_{i,\epsilon_i}(x)}.\]

    Thus, at each point $x$, we have
    \begin{eqnarray*}
        \Delta_f\rtilde &=& \sum_i\left[(\Delta_f \psi_i) r_{i,\epsilon_i} + \psi_i\Delta_f r_{i,\epsilon_i} + 2\inner{\grad \psi, \grad r_{i,\epsilon_i}}\right]\\
                    &=& \sum_i(\Delta_f \psi_i) r_{i,\epsilon_i} + \sum_i\psi_i\Delta_f r_{i,\epsilon_i} + b,\\
    \end{eqnarray*}
    and, since $\sum_i (\Delta_f \psi_i) r = r \Delta_f \left(\sum_i \psi_i\right) = r \Delta_f 1 = 0$,

    \[\Delta_f\rtilde = \sum_i(\Delta_f \psi_i) (r_{i,\epsilon_i} - r) + \sum_i\psi_i\Delta_f r_{i,\epsilon_i} + b.\]

    Let us now estimate the right side of above equation.

    By (\ref{controle_r_menos_r_i}), we have, for $x \in U_i$,
    \[|\Delta_f \psi_i (x) (r_{i,\epsilon_i}(x) - r(x))| \leq m_i \frac{\delta(x)}{2^im_i} \leq \frac{\delta(r(x))}{2},\]
    and, by Lemma \ref{controle_lap_f_r_teorema1_i} below,
    \[\Delta_f r_{i,\epsilon_i}(x) \leq \sup_{y \in B_x(1)}\{\Delta_f r(y)\} + \frac{\delta(r(x))}{2},\]
    so,
    \begin{eqnarray*}
        \sum_i\psi_i\Delta_f r_{i,\epsilon_i} &\leq& \sum_i \psi_i \left(\sup_{y \in B_x(1)}\{\Delta_f r(y)\}  + \frac{\delta(r(x))}{2}\right)\\
                                              &=&   \sup_{y \in B_x(1)}\{\Delta_f r(y)\} + \frac{\delta(r(x))}{2}.
    \end{eqnarray*}

    Therefore, for all $x \in M$,
    \[\Delta_f\rtilde (x) = \sup_{y \in B_x(1)}\{\Delta_f r(y)\} + \delta(x) + b(x).\]

    Since $\sum_i\inner{\grad \psi, \grad r} = \inner{\grad \left(\sum_i \psi_i\right), \grad r} = 0$ almost everywhere, we have
    \[b = 2\sum_i\inner{\grad \psi_i, (\grad r_{i,\epsilon_i} - \grad r)},\]
    almost everywhere. Thus, by (\ref{controle_grad_r_menos_grad_r_i}),

    \begin{eqnarray*}
        \|b\|_{L^1_f(M\setminus B_p(r))}
            &\leq& 2\|\sum_i\inner{\grad \psi_i, (\grad r_{i,\epsilon_i} - \grad r)}\|_{L^1_f(M\setminus B_p(r))}\\
            &\leq& 2\sum_{i\;|\;\,U_i \cap (M\setminus B_p(r)) \neq \emptyset}
                \left[\|\grad \psi_i\|_{L^\infty_f(U_i)}  \|\grad r_{i,\epsilon_i} - \grad r\|_{L^1_f(U_i)} \right] \\
            &\leq& 2\sum_{i\;|\;U_i \cap (M\setminus B_p(r)) \neq \emptyset}
                m_i\frac{\delta(r-1)}{2^im_i} \\
            &\leq& \delta(r-1).
    \end{eqnarray*}

    \begin{lemma}\label{controle_lap_f_r_teorema1_i}
        In the above construction, $\epsilon_i$ may be chosen so that
        \[\Delta_f r_{i,\epsilon_i}(x) \leq \sup_{y \in B_x(1)}\{\Delta_f r(y)\} + \frac{\delta(r(x))}{2},\]
        for all $x \in U_i$.
    \end{lemma}
    \proof
        Initially, notice that, by (\ref{controle_grad_r_i_pontual}) and by the second condition assumed about the neighborhoods $U_i$
        \begin{eqnarray*}
            |\inner{\grad f(x), \grad r(x)} - \langle\sum_j f_j\partial \textbf{x}^j(x), \grad r(x)\rangle |
                &=&
                    |\langle\grad f(x) - \sum_j f_j\partial \textbf{x}^j(x), \grad r(x)\rangle| \\
                &\leq&
                    \frac{\delta(r(x))}{2}
        \end{eqnarray*}
        and
        \begin{multline*}
            \displaystyle |\innerS{\grad f(x), \grad r_{i,\varepsilon_i}(x)} - \innerS{\sum_j f_j\partial \textbf{x}^j(x), \grad r_{i,\varepsilon_i}(x)}| = \\
                \begin{array}{cl}
                    \displaystyle
                        =& \displaystyle
                        |\langle\grad f(x) - \sum_j f_j\partial \textbf{x}^j(x), \grad r_{i,\varepsilon_i}(x)\rangle| \\
                        \leq&\displaystyle
                        \frac{\delta(r(x))}{2}.
                \end{array}
        \end{multline*}

        On the other hand, we have
        \begin{eqnarray*}
                \Delta_f r_{i,\varepsilon_i} (x)
            &=&
                \Delta r_{i,\varepsilon_i}(x) - \inner{\grad f(x), \grad r_{i,\varepsilon_i}(x)}\\
            &\leq&
                \Delta r_{i,\varepsilon_i}(x) - \inner{f_i \partial \textbf{x}^i(x), \grad r_{i,\varepsilon_i}(x)} + \frac{\delta(r(x))}{2}\\
            &=&
                \frac{\delta(r(x))}{2}\\
            &&  + \frac{1}{\epsilon^n}\int_{\R^n} \xi\left(\frac{\textbf{z}_i}{\epsilon}\right)
                \left[\Delta r(\textbf{y}_i + \textbf{z}_i) - \inner{f_i \partial \textbf{x}^i(\textbf{y}_i+\textbf{z}_i), \grad r(\textbf{y}_i+\textbf{z}_i) } \right]d\textbf{z}_i\\
            &=&
                \frac{\delta(r(x))}{2}\\
            &&  + \frac{1}{\epsilon^n}\int_{\R^n} \xi\left(\frac{\textbf{z}_i}{\epsilon}\right)
                \left[\Delta r(\textbf{y}_i+\textbf{z}_i) - \inner{\grad f(\textbf{y}_i+\textbf{z}_i), \grad r(\textbf{y}_i+\textbf{z}_i)} \right.\\
            &&  \phantom{\frac{1}{\epsilon^n}\int_{\R^n} \xi\left(\frac{\textbf{z}_i}{\epsilon}\right)[}
                + \left.\frac{\delta(r(\textbf{y}_i+\textbf{z}_i))}{2} \right]d\textbf{z}_i\\
            &=&
                \delta(r(x))\\
            &&  +\frac{1}{\epsilon^n}\int_{\R^n} \xi\left(\frac{\textbf{z}_i}{\epsilon}\right)
                \left[\Delta r(\textbf{y}_i+\textbf{z}_i) - \inner{\grad f(\textbf{y}_i+\textbf{z}_i), \grad r(\textbf{y}_i+\textbf{z}_i)} \right]d\textbf{z}_i\\
            &=&
                \delta(r(x))+ \frac{1}{\epsilon^n}\int_{\R^n} \xi\left(\frac{\textbf{z}_i}{\epsilon}\right)
                \Delta_f r(\textbf{y}_i+\textbf{z}_i)d\textbf{z}_i\\
            &\leq&
                \delta(r(x))+ \left(\sup_{y \in B_x(1)}\Delta_f r\right) \frac{1}{\epsilon^n}\int_{\R^n} \xi\left(\frac{\textbf{y}_i - \textbf{z}_i}{\epsilon}\right)d\textbf{z}_i\\
            &\leq&
                \delta(r(x))+ \sup_{y \in B_x(1)}\{\Delta_f r(y)\}.
        \end{eqnarray*}

    \proofend

    Define the lower level set $\widetilde{B}_p(t)$ of the function $\rtilde$ by
    \[\widetilde{B}_p(t) := \{x \in M, \tilde{r}(x) < t\}.\]

    We have the following relationship between the volumes of $B_p(t)$ and $\widetilde{B}_p(t)$.

    \begin{lemma}\label{comparacao_v_r_e_v_r_til}
        Given $a > 1$, the function $\delta$ in the construction of $\rtilde$ may be taken so that,
        \begin{enumerate}
            \item if $\vol_f(M) = \infty$,
                \[a^{-1} \vol_f(B_p(t)) \leq \vol_f(\widetilde{B}_p(t)) \leq a\vol_f(B_p(t));\]
            \item if $\vol_f(M) < \infty$,
                \begin{multline*}
                    a^{-1} (\vol_f(M) - \vol_f(\widetilde{B}_p(t))) \leq\\
                    \leq \vol_f(M) - \vol_f(B_p(t)) \leq \\
                    \leq a(\vol_f(M) - \vol_f(B_p(t)))
                \end{multline*}
        \end{enumerate}
            for all $t > 1$.
    \end{lemma}

    \proof
        We will proof first the case $\vol_f(M) = \infty$.

        For any $t > 1$, take $\eta_t$ such that
        \[a^{-1}\vol_f(B_p(t + \eta_t)) < \vol_f(B_p(t)) < a\vol_f(B_p(t - \eta_t)).\]

        By continuity, we suppose that the above inequalities hold in an open neighborhood $U_t$ of $t$. By considering all $t \geq 1$, we have a cover $\{U_t\}$ of $[1,\infty)$, from which we extract a locally finite cover $\{U_{t_i}\}$. Now, take the function $\delta$ so that $2\delta(t) \leq \inf_{\{i, t \in U_{t_i}\}} \eta_{t_i}$, $\delta(t) < 1$ and $\delta(t-1) < 2\delta(t)$.

        If $x \in \widetilde{B}_p(t)$, we have $\rtilde(x) < t$, and so $r(x) < t + \delta(t)$. Indeed, if \linebreak this was not the case, we would have $r(x) \geq t + \delta(t)$ and so \linebreak $|r(x) - \rtilde(x)| \geq r(x) - \rtilde(x) \geq t + \delta(t) - \rtilde(x) > \delta(t)$, contradicting the construction of $\rtilde$. Therefore $\widetilde{B}_p(t) \subset B_p(t + \delta)$.

        Now take $x \in B_p(t - 2\delta(t))$. If $r(x) < t - 1$ then, by the construction of $\rtilde$, we have $\rtilde(x) - r(x) < \delta(r(x)) < 1$ and so $\rtilde(x) < 1 + r(x) <\linebreak 1 + t -1 = t$. On the other hand, if $t-1 \leq r(x) < t - 2\delta(t)$ then \linebreak $\rtilde(x) - r(x) <  \delta(r(x)) < \delta(t-1) < 2\delta(t)$ and so $\rtilde(x) < r(x) + 2\delta(x) < t$. Therefore $B_p(t - 2\delta(t)) \subset \widetilde{B}_p(t)$.

        We have just proved that
        \[B_p(t - \eta_{t_i}) \subset B_p(t - 2\delta(t)) \subset \widetilde{B}_p(t) \subset B_p(t + 2\delta(t)) \subset B_p(t + \eta_{t_i}),\]
        thus
        \begin{eqnarray*}
        a^{-1}\vol_f(B_p(t)) &\leq& \vol_f(B_p(t - \eta_{t_i}))\\
                             &\leq& \vol_f(\widetilde{B}_p(t))\\
                             &\leq& \vol_f(B_p(t + \eta_{t_i}))\\
                             &\leq& a\vol_f(B_p(t)).
        \end{eqnarray*}

        For the second case ($\vol_f(M) = \infty$), if we take $\eta_t$ so that
        \begin{eqnarray*}
            a^{-1}(\vol_f(M) - \vol_f(B_p(t - \eta_t)))
                &<& \vol_f(M) - \vol_f(B_p(t))\\
                &<& a(\vol_f(M) - \vol_f(B_p(t + \eta_t))),
        \end{eqnarray*}
        the result follows as above.
    \proofend

    Since it will be used later, we present an analogous construction for $L^1$ and $\Delta$. Given $\delta: \R^+ \rightarrow \R^+$, we construct $b$ and $\rtilde$ in $C^\infty(M)$ such that
    \begin{enumerate}[label={\roman*.}]
        \item $\displaystyle \norm{b}_{L^1(M \setminus B_p(r))} \leq \delta(r)$;
        \item $\displaystyle \norm{\grad\rtilde - \grad r}_{L^1(M \setminus B_r(p))} \leq \delta(r)$;
        \item $\displaystyle \module{\rtilde(x) - r(x)} \leq \delta(r(x))$ and $\module{\grad\rtilde(x)}\leq 2, \forall x \in M, r(x) > 2$, and
        \item $\displaystyle \Delta\rtilde(x) \leq \sup_{y \in B_x(1)}\{\Delta r(y)\} + \delta(r(x)) + \module{b(x)}, \forall x \in M, r(x) > 2$  in the sense of distributions.
    \end{enumerate}

    Using this new construction of $\rtilde$, Lemma \ref{comparacao_v_r_e_v_r_til} remais valid if we replace $\vol_f$ with $\vol$.

\section{Estimates for $\Delta_f r$ and $\vol_f(B_p(r))$ in the hypotheses of Theorem \ref{principal1ext}}\label{estimativas_vol_f_teorema1ext}

    Before we can prove Theorem \ref{principal1ext}, we need some results concerning the \linebreak $f$-volume of $B_p(r)$. In this section, we obtain estimates for the growth of the function $r \mapsto \vol_f(B_p(r))$, in the case $\vol_f(M) = \infty$, and its decay, in the case $\vol_f(M) < \infty$.

    For the sake of simplicity, in the proofs of this section, we define the notation $V_f(r):=\vol_f(B_p(r))$.

    \begin{lemma}\label{controle_de_lap_r_teorema1ext}
        In the hypothesis of Theorem \ref{principal1ext}, if $\ric_f(\grad r, \grad r) \geq 0$ and \linebreak $\lim_{r \rightarrow \infty} \frac{|f|}{r} = 0$, then
        \[\lim_{r\rightarrow \infty} \Delta_f r \leq 0.\]
    \end{lemma}

    \proof
        Observe that
        \[0 \leq \ric_f (\grad r,\grad r) = \ric(\grad r,\grad r) + \grad_{\grad r}\grad_{\grad r} f,\]
        so
        \begin{equation}\label{controle_de_lap_r_teorema1ext_c1}
            - \grad_{\grad r}\grad_{\grad r} f \leq \ric(\grad r,\grad r).
        \end{equation}

        Let $\gamma(t)$ be a minimizing normal geodesic starting from $p$. By the Bochner Formula, we have, along $\gamma$,
        \[0 = \frac{1}{2}\Delta|\grad r| = |\hess r|^2 + (\Delta r)' + \ric(\grad r,\grad r),\]
        where the derivatives are calculated with respect to $t$. Hence
        \begin{equation}\label{controle_de_lap_r_teorema1ext_c2}
           (\Delta r)' = - |\hess r|^2 - \ric(\grad r,\grad r) \leq - \frac{(\Delta r)^2}{n-1} - \ric(\grad r,\grad r).
        \end{equation}

        Combining (\ref{controle_de_lap_r_teorema1ext_c1}) and (\ref{controle_de_lap_r_teorema1ext_c2}), we have
        \begin{eqnarray*}
            (r^2\Delta_f r)'
            &=& (r^2\Delta r - r^2\grad_{\grad r} f)'\\
            &=& 2r\Delta r + r^2(\Delta r)' - 2r\grad_{\grad r} f - r^2\grad_{\grad r}\grad_{\grad r} f\\
            &\leq& 2r\Delta r - r^2\;\frac{(\Delta r)^2}{n-1} - r^2\ric(\grad r, \grad r) - 2r\grad_{\grad r} f + r^2 \ric(\grad r, \grad r)\\
            &=& 2r\Delta r - r^2\;\frac{(\Delta r)^2}{n-1} - 2r\grad_{\grad r} f\\
            &=& -\left(\frac{r\Delta r}{\sqrt{n-1}} - \sqrt{n-1}\right)^2 + (n-1) - 2r\grad_{\grad r} f\\
            &\leq& (n-1) - 2r\grad_{\grad r} f.
        \end{eqnarray*}

        Since $r^2\Delta_f r$ vanishes for $r = 0$, we have
        \begin{eqnarray*}
            r^2\Delta_f r &\leq& \int_0^r (n-1)\;dt - 2\int_0^r t\grad_{\grad r} f\;dt\\
            &=& (n-1)r + 2\int_0^r f(t) \;dt - 2rf(r).
        \end{eqnarray*}

        Taking $r_0$ such that $|f| < \epsilon r$ for all $r > r_0$, we get
        \begin{eqnarray*}
            r^2\Delta_f r &<& (n-1)r + 2\epsilon \int_{r_0}^r t\;dt + 2\int_0^{r_0} f(t)\;dt + 2\epsilon r^2.\\
                       &\leq& (n-1)r + \epsilon r^2  - \epsilon r_0^2 + r_0 C_{r_0} + 2\epsilon r^2,
        \end{eqnarray*}
        where $C_{r_0} = \sup_{B_p(r_0)} |f|$. Therefore,
        \[\Delta_f r < \frac{(n-1)}{r} + 3\epsilon + \frac{r_0 C_{r_0}}{r^2}.\]

        Since the above expression is independent on the geodesic we consider, we can take $r_1 > 0$ such that $\Delta_f r < 4\epsilon$  for $x \in M\setminus B_p(r_1)$. Therefore, $\lim_{r \rightarrow \infty} \Delta_f r  \leq 0$.
    \proofend

    \begin{lemma}\label{estimativa_vol_infty_teorema1}
        In the hypotheses of Theorem \ref{principal1ext} and supposing $\vol_f(M) = \infty$, for all $\epsilon > 0$ and $t_0$ large enough, there exists $R = R(\epsilon, t_0)$ such that, for $t > R$,
        \[\int_{B_p(t)\setminus B_p(t_0)} \module{\Delta_f \rtilde}\;e^{-f}d\vol_g \leq \epsilon \vol_f(B_p(t + 1)).\]
    \end{lemma}

    \proof
        By the construction we have done in section \ref{rtil_teorema1},
        \[\Delta_f\rtilde \leq \sup_{y \in B_x(1)}\{\Delta_f r(y)\} + \delta(r(x)) + b(x),\]
        so,
        \[|\Delta_f\rtilde| \leq 2 \left(\sup_{y \in B_x(1)}\{\Delta_f r(y)\} + \delta(r(x)) + |b(x)|\right) - \Delta_f\rtilde,\]
        in the sense of distributions.

        Given $\epsilon > 0$, by Lemma \ref{controle_de_lap_r_teorema1ext} we take $r_0 > 0$ such that, for $r(x) \geq r_0$, $\sup_{y \in B_x(1)}\{\Delta_f r(y)\} + \delta(r(x)) < \epsilon$. Thus, for $r > r_0$,
        \begin{eqnarray*}
            \int_{B_p(r)\setminus B_p(r_0)} |\Delta_f \rtilde|\;e^{-f}d\vol_g
            &\leq& \frac{\epsilon}{2}(V_f(r) - V_f(r_0) + 2\int_{M\setminus B_p(r_0)} |b|\;e^{-f}d\vol_g\\
            && - \int_{B_p(r) \setminus B_p(r_0)} \Delta_f \rtilde\;e^{-f}d\vol_g\\
            &=& \frac{\epsilon}{2}(V_f(r) - V_f(r_0)) + 2\int_{M\setminus B_p(r_0)} |b|\;e^{-f}d\vol_g\\
            &&  - \int_{\partial B_p(r)}\frac{\partial \rtilde}{\partial n}\;e^{-f}d\vol_g + \int_{\partial B_p(r_0)}\frac{\partial \rtilde}{\partial n} \;e^{-f}d\vol_g\\
        \end{eqnarray*}

        Since $-\alpha \leq |\alpha - 1|, \forall \alpha \in \R$, we get
        \[- \int_{\partial B_p(r)}\frac{\partial \rtilde}{\partial n}\;e^{-f}d\vol_g \leq \int_{\partial B_p(r)}\left|\frac{\partial \rtilde}{\partial n}-1\right|\;e^{-f}d\vol_g,\]
        therefore,
        \begin{eqnarray}
            \nonumber
                \int_{B_p(r)\setminus B_p(r_0)} |\Delta_f \rtilde|\;e^{-f}d\vol_g
                    &\leq&
                        \frac{\epsilon}{2}(V_f(r) - V_f(r_0))\\
            \nonumber
                        && + 2\int_{M\setminus B_p(r_0)} |b|\;e^{-f}d\vol_g\\
            \nonumber
                        && + \int_{\partial B_p(r_0)}\frac{\partial \rtilde}{\partial n}\;e^{-f}d\vol_g\\
            \label{parcial_estimativa_vol_infty_teorema1}
                        && + \int_{\partial B_p(r)}\left|\frac{\partial \rtilde}{\partial n}-1\right|\;e^{-f}d\vol_g.
        \end{eqnarray}

        Since
        \[\int_{M\setminus B_p(r_0)} |b| \;e^{-f}d\vol_g \leq \delta(r_0) \leq \epsilon,\]
        and $\vol_f(M) = \infty$, fixing $r_0$ such that
        \begin{gather}
            \label{vol_f_r0_maior_que_2}
            V_f(r_0) \geq 1,\\
            \label{grad_rtil_menos_grad_r}
            \|\grad \rtilde - \grad r\|_{L^1_f(M\setminus B_p(r_0))} \leq 1.
        \end{gather}

        Taking $r$ large enough, we have
        \[\epsilon (V_f(r) - V_f(r_0)) \geq \epsilon + \int_{\partial B_p(r_0)}\frac{\partial \rtilde}{\partial n}\;e^{-f}d\vol_g,\]
        thus
        \begin{eqnarray*}
            \int_{B_p(r)\setminus B_p(r_0)} |\Delta_f \rtilde|\;e^{-f}d\vol_g
                &\leq&
            \frac{\epsilon}{2}(V_f(r) - V_f(r_0)) + \epsilon + \int_{\partial B_p(r_0)}\frac{\partial \rtilde}{\partial n}\;e^{-f}d\vol_g\\
             && + \int_{\partial B_p(r)}\left|\frac{\partial \rtilde}{\partial n}-1\right|\;e^{-f}d\vol_g \\
                &\leq&
            \frac{3\epsilon}{2}V_f(r) + \int_{\partial B_p(r)}\left|\frac{\partial \rtilde}{\partial n}-1\right|\;e^{-f}d\vol_g.\\
        \end{eqnarray*}

        By (\ref{grad_rtil_menos_grad_r}),
        \begin{eqnarray*}
            \int_{r}^{r+1} dt \int_{\partial B_p(t)} \left|\frac{\partial \rtilde}{\partial n}-1\right|\;e^{-f}d\vol_g
            &=&
            \int_{r}^{r+1} dt \int_{\partial B_p(t)} \left|\inner{\grad \rtilde - \grad r, \grad r}\right|\;e^{-f}d\vol_g\\
            &\leq&
            \int_{r}^{r+1} dt \int_{\partial B_p(t)} \left|\grad \rtilde - \grad r\right|\;e^{-f}d\vol_g\\
            &\leq&
            \int_{M\setminus\partial B_p(r)} \left|\grad \rtilde - \grad r\right|\;e^{-f}d\vol_g\\
            &\leq&
            \int_{M\setminus\partial B_p(r_0)} \left|\grad \rtilde - \grad r\right|\;e^{-f}d\vol_g\\
            &\leq&
            1,
        \end{eqnarray*}
        therefore, for some $r' \in [r, r+1]$, we have
        \[\int_{\partial B_p(t)} \left|\frac{\partial \rtilde}{\partial n}-1\right|\;e^{-f}d\vol_g \leq 1.\]

        Hence,
        \begin{eqnarray*}
            \int_{B_p(r)\setminus B_p(r_0)} |\Delta_f \rtilde|\;e^{-f}d\vol_g
            &\leq&
            \int_{B_p(r')\setminus B_p(r_0)} |\Delta_f \rtilde|\;e^{-f}d\vol_g\\
            &\leq&
            \frac{3\epsilon}{2}(V_f(r') - V_f(r_0) - \int_{\partial B_p(r')}\frac{\partial \rtilde}{\partial n}\;e^{-f}d\vol_g\\
            &\leq&
            \frac{3\epsilon}{2}(V_f(r') - V_f(r_0) + 1\\
            &=&
            \frac{3\epsilon}{2}V_f(r') - \frac{3\epsilon}{2}V_f(r_0) + 1\\
            &\leq&
            \frac{3\epsilon}{2}V_f(r')\\
            &\leq&
            \frac{3\epsilon}{2}\vol_f(B_p(r+1)).
        \end{eqnarray*}

        We have used (\ref{vol_f_r0_maior_que_2}) in the last second line.

    \proofend

    The following result estimates the growth of $\vol_f(B_p(r))$.

    \setcounter{theorem}{2}
    \begin{theorem}
        Let $M$ be a non-compact complete manifold. If $f: M \rightarrow \R$ is a smooth function such that $\ric_f \geq 0$ and $\lim_{r \rightarrow \infty}\frac{|f|}{r} = 0$, then
        \[\lim_{r \rightarrow \infty}\frac{\vol_f(B_p(r))}{e^{Cr}} = 0,\]
        for all $C > 0$.
    \end{theorem}

    \proof
        Denoting $A_f(r) = \vol_f(\partial B_p(r))$, by Lemma \ref{controle_de_lap_r_teorema1ext},

        \begin{eqnarray*}
            A_f'(r) &=& \dd{}{r} \int_{\partial B_p(r)} e^{-f} d\vol_g\\
                    &=& - \int_{\partial B_p(r)} \inner{\grad f, \grad r}e^{-f} d\vol_g + \int_{\partial B_p(r)} \Delta r\;e^{-f} d\vol_g \\
                    &=& \int_{\partial B_p(r)} \Delta_f r\;e^{-f} d\vol_g\\
                    &\leq& h(r) A_f(r),
        \end{eqnarray*}
        where $h(r) := \sup_{\partial B_p(r)}\Delta_f r$. Thus,
        \[(\log A_f(r))' = \frac{A_f'(r)}{A_f(r)} \leq h(r),\]
        which implies
        \[\log A_f(r) \leq  \int_0^r h(t)\;dt + \log A_f(1),\]
        and therefore,
        \[A_f(r) \leq C_1 e^{\int_0^r h(t)\;dt}.\]

        Given $C > 0$, if $r_0 > 0$ is such that $h(r) < C/2$ for $r \geq r_0$, we have
        \[A_f(r) < C_1 e^{\int_0^{r_0} h(t)\;dt + \frac{C}{2}(r - r_0)} = C_2 e^{\frac{C}{2}r},\]
        for all $r > r_0$. Therefore,
        \begin{eqnarray*}
            \lim_{r \rightarrow \infty} \frac{\vol_f(B_p(r))}{e^{C r}}
                &=& \lim_{r \rightarrow \infty} \frac{A_f(r)}{C e^{C r}}\\
                &<& \lim_{r \rightarrow \infty} \frac{C_2 e^{\frac{C}{2}r}}{C e^{C r}}\\
                &=& \frac{C_2}{C}\lim_{r \rightarrow \infty} e^{-\frac{C}{2}r} = 0.\\
        \end{eqnarray*}

    \proofend

     \begin{lemma}\label{controle_crescimento_volume}
        Let $f$ (not necessarily as in the hypotheses of Theorem \ref{principal1ext}) be such that \[\lim_{r \rightarrow \infty}\frac{\vol_f(B_p(r)}{e^{Cr}} = 0,\]
        for all $C \in \R$. Then, for $R > 0$, there exists a sequence $r_k \rightarrow \infty$ of real numbers such that $\vol_f(B_p(r_k + R + 1)) \leq 2\vol_f(B_p(r_k))$.
    \end{lemma}

    \proof
        Given $R > 0$, set $r_k = (k-1)(R+1)$, $k \in \N$. If there exists infinitely many indices $k$ satisfying $V_f(r_k + R + 1) \leq 2V_f(r_k)$,
        by passing to a subsequence of $r_k$, we will have the desired sequence. If this is not the case, we have
        \[V_f(k(R + 1)) = V_f(r_k + R + 1) > 2V_f(r_k) = 2V_f((k-1)(R+1)),\]
        for $k$ larger than some $k_0$. Therefore,
        \[V_f(k(R+1)) > 2^{k - k_0}V_f((k_0-1)(R+1)).\]

        Taking $C$ such that $2^{k - k_0}V_f((k_0-1)(R+1)) > e^{Ck(R+1)}$ for large $k \in \N$, we have
        \[\frac{V_f(k(R+1))}{e^{Ck(R+1)}} > 1,\]
        which contradicts the hypotheses.
    \proofend

    The last two Lemmas lead us to the following.

    \begin{lemma}\label{controle_crescimento_volume_teorema1ext}
        In the hypotheses of Theorem \ref{principal1ext}, for $R > 0$, there exists a sequence $r_k \rightarrow \infty$ of real numbers such that $\vol_f(B_p(r_k + R + 1)) \leq 2\vol_f(B_p(r_k))$.
     \end{lemma}

    In the hypothesis of Theorem \ref{principal1ext} we cannot prove that the $f$-volume of $M$ is infinite. Supposing now $\vol_f(M) < \infty$, we will get an interesting estimate on the decay of the $f$-volume. \\

     If $\vol_f(M) < \infty$, we say that the $f$-volume of $M$ decays exponentially at $p$ if exists $C > 0$ such that
     \[\vol_f(M) - \vol_f(B_p(r)) < e^{-Cr},\]
     for all $r$ large enough.

     \begin{theorem}
        Let $M$ be a non-compact complete manifold. If $f: M \rightarrow \R$ is a smooth function such that $\ric_f \geq 0$, $\lim_{r \rightarrow \infty}\frac{|f|}{r} = 0$ and $\vol_f(M) < \infty$, then the $f$-volume of $M$ does not decay exponentially at $p$.
     \end{theorem}

     \proof
        This proof is inspired by the idea of \cite[Lemma 5.2]{munteanu}.

        Suppose in contradiction the existence of $C > 0$ such that
        \begin{equation}\label{absurdo_f_volume_decai_exp}
            \vol_f(M \setminus B_p(r)) \leq e^{-Cr}.
        \end{equation}

        Take $R > 0$ such that $\frac{|f|}{r} < \epsilon < C/4$ in $M\setminus B_p(R-1)$, and let $q \in M$ be such that $d(q,p) = R + S$, for arbitrary $S > 0$. Denote by
        \[J_f(t,\xi)\;dt\,d\xi := e^{-f}d\vol_g|_{\exp_q(t\xi)}\]
        the $f$-volume form in geodesic polar coordinates $(t, \xi)$ from $q$. We have
        \[d\vol_g = e^f\;J_f(t,\xi)\;dt\,d\xi,\]
        so
        \begin{eqnarray}
            \label{d_forma_de_fvolume}
            \Delta t &=& \dd{}{t} \log(e^f\;J_f(t,\xi))\\
            \nonumber &=& \dd{}{t} (f+\log(J_f(t,\xi) )\\
            \nonumber &=& \inner{\grad f, \grad t} + \dd{}{t} \log(J_f(t,\xi)).
        \end{eqnarray}

        Let $p_R$ be a point in the minimizing geodesic from $p$ to $q$, and such that $d(p_R, p) = R$. Taking $a$, $-1 < a < 1$, let $x$ be an arbitrary point of $\partial B_{q}(S+a) \cap B_{p_R}(1)$ and $\gamma$ be the minimizing geodesic such that $\gamma(0) = q$, $\gamma(S+a) = x$. Notice that, since $\gamma$ is a minimizing geodesic, if we define $y:=\gamma(2+a)$, each $x$ will determine, in an injective way, $y \in B_q(2+a)$. Moreover, we have $y = \exp_q(a_y\,\xi_y)$, where $a_y = a$ and $\xi_y = \gamma'(0)$.

        \begin{center}\includegraphics[height=5cm]{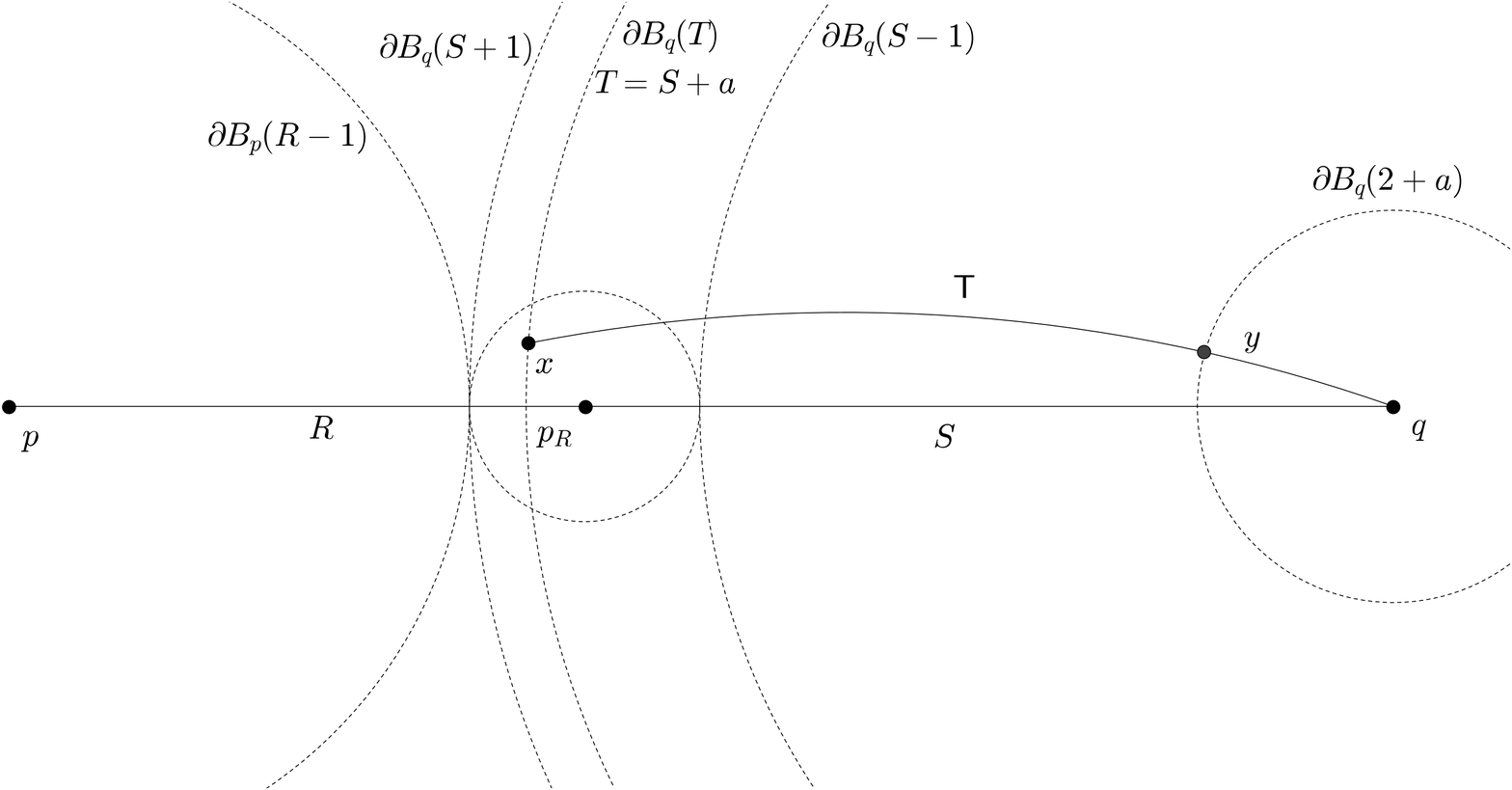}\end{center}

        For the sake of simplicity, we denote $T = S+a$ and $f(t) = f(\gamma(t))$. Since
        \[r(\gamma(t)) = d(p, \gamma(t)) \leq d(p, x) + d(x, \gamma(t)) \leq R + 1 + T - t,\]
        we have
        \begin{equation}\label{estimativa_f_gamma_t_teo1ext}
            |f(\gamma(t))| < \epsilon r(\gamma(t)) \leq \epsilon((R + 1 + T) - t).
        \end{equation}

        Since $\ric + \grad \grad f = \ric_f \ge 0$, we have
        \begin{eqnarray*}
            0 = \Delta \frac{|\grad t|^2}{2} &=& |\hess t|^2 + (\Delta t)' + \ric(\grad t, \grad t)\\
                &\geq& \frac{(\Delta t)^2}{n-1} + (\Delta t)' - f'',\\
        \end{eqnarray*}
        so
        \[(\Delta t)' \leq f'' - \frac{(\Delta t)^2}{n-1}.\]

        Multiplying by $t^2$ and integrating from $0$ to $s \leq T$, we get
        \[\int_0^s t^2(\Delta t)' \;dt \leq \int_0^s t^2 f''(t)\;dt - \frac{1}{n-1}\int_0^s t^2(\Delta t)^2\;dt,\]
        which is equivalent to
        \begin{eqnarray*}
             s^2 \Delta t (s)
                &\leq& - \frac{1}{n-1}\int_0^s \left[ t^2(\Delta t)^2 - 2(n-1)t \Delta t \right]\;dt + \int_0^s t^2 f''(t)\;dt\\
                &=& - \frac{1}{n-1}\int_0^s \left[t\Delta t - (n-1)\right]^2\;dt + (n-1)\int_0^s dt + \int_0^s t^2 f''(t)\;dt\\
                &\leq&(n-1)s + s^2 f'(s) - 2\int_0^s t f'(t)\;dt.
        \end{eqnarray*}

        Therefore,
        \[\Delta t (s) \leq \frac{n-1}{s} + f'(s) - \frac{2}{s^2}\int_0^s t f'(t)\;dt.\]

        Using (\ref{d_forma_de_fvolume}) and integrating again from $b:=2+a$ to $T$,
        \begin{eqnarray}
         \nonumber   \log(J_f(T,\xi)) - \log(J_f(b,\xi)) &=& \int_b^T \dd{}{s} \log(J_f(s, \xi))\;ds\\
         \nonumber   &=& \int_b^T (\Delta t(s) - f'(s))ds\\
         \nonumber   &=& (n-1)(\log T - \log b)\\
         \label{d_forma_de_fvolume1} && - \int_b^T \left(\frac{2}{s^2}\int_0^s t f'(t) dt\right)ds.
        \end{eqnarray}

        Now we estimate the last term in (\ref{d_forma_de_fvolume1}). Integrating by parts and \linebreak using (\ref{estimativa_f_gamma_t_teo1ext}),
        \begin{eqnarray*}
            - \int_b^T \frac{2}{s^2}ds\int_0^s t f'(t)\;dt
                &=& \left.\left(\frac{2}{s}\int_0^s t f'(t)\;dt\right)\right|_b^T - 2(f(T) - f(b))\\
                &=& - \frac{2}{T}\int_0^T f(t)\;dt + \frac{2}{b}\int_0^b f(t)\;dt\\
                &\leq& \frac{2}{T}\int_0^T |f(t)|\;dt + \frac{2}{b}\int_0^b |f(t)|\;dt\\
                &\leq& 2\epsilon\left[\frac{1}{T}\int_0^T ((R+1+T) - t)\;dt\right.\\
                && + \left.\frac{1}{b}\int_0^b ((R+1+T) - t)\;dt\right]\\
                &=& 2\epsilon\left[2(R+1+T) - \frac{T + b}{2}\right]\\
                &\leq& 4\epsilon(R+1) + 4\epsilon T.\\
        \end{eqnarray*}

        Putting this into (\ref{d_forma_de_fvolume1}), we get
        \[
            \log(J_f(T,\xi)) - \log(J_f(b,\xi)) \leq
                (n-1)\log T + 4\epsilon T + 4\epsilon(R+1).
                \]

        Reminding that $S - 1 \leq T \leq S + 1$ and $-1 \leq a \leq 1$, there exists $c_1$ and $c_2$ constants and independent from $a, \xi, S, q$ such that
        \begin{eqnarray}\label{forma_de_fvolume}
            J_f(a+2,\xi) &=& J_f(b, \xi) \\
            \nonumber &\geq& e^{- (n-1)\log T - 4\epsilon T - 4\epsilon(R+1)}J_f(T,\xi) \\
            \nonumber &\geq& c_1 e^{-c_2\log S - 4\epsilon S} J_f(S+a,\xi).
        \end{eqnarray}

        Now, let $U\subset(B_q(3))$ and $V \subset T_qM$ be, respectively, the set of all $y$ and $(a_y, \xi_y)$ constructed as above. Notice that each $x \in B_q(1)$ determines an unique $y$ in $U$ and $(a_y, \xi_y)$ in $V$. Therefore
        \begin{eqnarray*}
            e^{-C(R+S-3)} &\geq&
                \vol_f(M \setminus B_p(R+S-3)) \\
                &\geq& \vol_f(B_q(3))\\
                &\geq& \vol_f(U)\\
                &\geq& \int_V J_f(a,\xi)\;dr\;d\xi\\
                &\geq& c_1 e^{-c_2\log S - 4\epsilon S} \int_V J_f(S+a,\xi)\;dr\;d\xi\\
                &\geq& c_1 e^{-c_2\log S - 4\epsilon S} \int_{exp_q(V)}\;e^{-f}d\vol_g\\
                &\geq& c_1 e^{-c_2\log S - 4\epsilon S} \vol_f(B_{p_R}(1)).
        \end{eqnarray*}

        Since $M$ is complete, we take $K = \min_{p_R} \vol_f(B_{p_R}(1))$ and so
        \begin{eqnarray*}
            \left(e^{-C(R-3)}\right)\;e^{-CS} &=& e^{-C(R+S-3)}\\
                                              &\geq&  K c_1 e^{-c_2\log S - 4\epsilon S}\\
                                              &\geq& \left(K c_1 e^{-c_2\log S}\right)\; e^{CS/2},
        \end{eqnarray*}
        which is an absurd for large $S$.

     \proofend

    \begin{lemma}\label{estimativa_vol_finito_teorema1ext}
        Suppose $\vol_f(M) < \infty$. Given $\epsilon > 0$, the construction of $\rtilde$ in Section \ref{rtil_teorema1} can be made so that, for $R > 0$ large enough, we have
        \[\int_{M_p\setminus B_p(t)} \module{\Delta_f \rtilde}\;e^{-f}d\vol_g \leq \epsilon (\vol_f(M) - \vol_f(B_p(t))) + 2\vol_f(\partial B_p(t)),\]
        for $t > R$.
    \end{lemma}

    \proof
        As remarked in the beginning of Section \ref{rtil_teorema1}, given $\varepsilon$ we choose $\delta$ so that
        \[\delta(r) \leq \epsilon(\vol_f(M) - V_f(r)).\]

        Therefore,
        \[\int_{M\setminus B(t)} |b| \leq \delta(t) \leq \epsilon(\vol_f(M) - V_f(t)).\]

        Since \[\left|\frac{\partial \rtilde}{\partial r} - 1\right| = |\inner{\grad \rtilde - \grad r, \grad r}|\leq |\grad \rtilde - \grad r|\]
        and
        \[\lim_{r\rightarrow +\infty}\|\grad \rtilde - \grad r\|_{L^1_f(M\setminus B_p(r))} = 0,\]
        by the hypotheses $\vol_f(M) < \infty$ we can take a sequence $r_i\rightarrow +\infty$ such that

        \[\lim_{i\rightarrow \infty} \int_{\partial B_p(r_i)} \left|\frac{\partial \rtilde}{\partial r} - 1\right|\;e^{-f}d\vol_g = 0\]
        and
        \[\lim_{i\rightarrow \infty} V_f(r_i) = \vol_f(M).\]

        Hence, taking $r_i \rightarrow \infty$ in inequality (\ref{parcial_estimativa_vol_infty_teorema1}) of the proof of Lemma \ref{estimativa_vol_infty_teorema1}, we have,
        \begin{eqnarray*}
                \int_{M\setminus B_p(r_0)} |\Delta_f \rtilde|\;e^{-f}d\vol_g
                    &\leq&
                \frac{\epsilon}{2}(\vol_f(M) - V_f(r_0)) + 2\int_{M\setminus B_p(r_0)} |b|\;e^{-f}d\vol_g\\
                && + \int_{\partial B_p(r_0)}\frac{\partial \rtilde}{\partial n}\;e^{-f}d\vol_g\\
                    &\leq&
                \frac{5\epsilon}{2}(\vol_f(M) - V_f(r_0))
                + \int_{\partial B_p(r_0)}\frac{\partial \rtilde}{\partial n}\;e^{-f}d\vol_g.
        \end{eqnarray*}

        Thereby, as $|\grad \rtilde| < 2$, we have
        \[\int_{M\setminus B_p(r_0)} |\Delta_f \rtilde|\;e^{-f}d\vol_g \leq \frac{5\epsilon}{2}(\vol_f(M) - V_f(r_0)) + 2\vol_f(\partial B_p(r_0)).\]

    \proofend

    \begin{lemma}\label{estimativa_vol_finito2_teorema1ext}
        Let $M$ and $f$ be as in the hypotheses of Theorem \ref{principal1ext} and suppose $\vol_f(M) < \infty$. Given $\epsilon > 0$, $C > 0$, there exist $R > 0$ large and a sequence of real numbers $r_k \rightarrow \infty$ such that
        \begin{multline*}
            \epsilon[\vol_f(M) - \vol_f(B_p(r_k - R))] + C\;\vol_f(\partial B(r_k - R)) \leq\\ 2\epsilon\;[\vol_f(M) - \vol_f(B(r_k))].
        \end{multline*}

    \end{lemma}

    \proof
        If there were not such $R$ and sequence $r_k$, for all $R>0$ and $r$ greater then some $r_0$ we would have
        \[ \epsilon(\vol_f(M) - V_f(r - R)) + C\;\vol_f(\partial B(r - R)) > 2\epsilon\;(\vol_f(M) - V_f(r)),\]
        so
        \begin{multline*}
            \epsilon\;e^{-\epsilon r/C}(\vol_f(M) - V_f(r - R)) + Ce^{-\epsilon r/C}\;\vol_f(\partial B(r - R)) >\\
            2\epsilon\;e^{-\epsilon r/C}\;(\vol_f(M) - V_f(r)),
        \end{multline*}
        and, since $[\vol_f(M) - V_f(r - R)]' = \vol_f(\partial B(r - R))$,
        \begin{equation}\label{parcial_lemma_estimativa_vol_finito_teorema1}
        \left[-e^{-\epsilon r/C}(\vol_f(M) - V_f(r - R))\right]^\prime > 2\frac{\epsilon}{C}\;e^{-\epsilon r/C}\;(\vol_f(M) - V_f(r)).
        \end{equation}

        Hence,
        \begin{multline*}
            \displaystyle \begin{array}{rc}
            \displaystyle e^{-\epsilon r/C}(\vol_f(M) - V_f(r - R))-\\
            \displaystyle - e^{-\epsilon (r+R)/C}(\vol_f(M) - V_f(r + R)) &=
            \end{array}\\
            \displaystyle \begin{array}{cl}
                =& \displaystyle \int_{r}^{r+R}\left[-e^{-\epsilon t/C}(\vol_f(M) - V_f(t - R))\right]^\prime dt\\
                >& \displaystyle 2\frac{\epsilon}{C}\int_{r}^{r+R}\;e^{-\epsilon t/C}\;(\vol_f(M) - V_f(t))dt\\
                >& \displaystyle 2\frac{\epsilon}{C}(\vol_f(M) - V_f(r+R))\int_{r}^{r+R}\;e^{-\epsilon t/C}dt\\
                =& \displaystyle -2(\vol_f(M) - V_f(r+R))[e^{-\epsilon (r+R)/C} - e^{-\epsilon r/C}],
            \end{array}
        \end{multline*}
        and so
        \[e^{-\epsilon r/C}(\vol_f(M) - V_f(r - R)) > (\vol_f(M) - V_f(r+R))[-e^{-\epsilon (r+R)/C} + 2e^{-\epsilon r/C}],\]
        which, after dividing by $e^{-\epsilon r/C}$, leads to
        \[(\vol_f(M) - V_f(r - R)) > (\vol_f(M) - V_f(r+R))[-e^{-\epsilon R/C} + 2].\]

        Taking $R$ large, we have $-e^{-\epsilon R/C} + 2 > L > 1$ and thus
        \[(\vol_f(M) - V_f(r - R)) > L (\vol_f(M) - V_f(r+R)),\]
        and, iterating,
        \[(\vol_f(M) - V_f(r - R)) > L^n (\vol_f(M) - V_f(r + (2n-1)R)).\]

        Therefore,
        \[\vol_f(M) - V_f(r + (2n-1)R) < L^{-n}(\vol_f(M) - V_f(r - R)) < L^{-n}(\vol_f(M)).\]

        The above inequality shows that the $f$-volume of $M$ decays exponentially, contradicting the Theorem \ref{decaimento_nao_exponencial_teorema1ext}.

    \proofend

\section{Proof of Theorem \ref{principal1ext}}\label{prova_do_teorema1ext}

    Now we have studied all the necessary tools, we can proceed with the proof of the Theorem \ref{principal1ext}.

    According to Corollary \ref{weyl_lap_f}, in order to prove that $\sigma_{\mathrm{ess}}(\Delta_f) = [0, \infty)$ in $L^2$, we need only to construct, for all $\lambda > 0$, a sequence $\{u_i\} \subset C_0^\infty(M)$ satisfying
    \begin{enumerate}
        \item $\displaystyle \frac{\lnorm{u_i}{\infty}{M}\;\lfnorm{(-\Delta_f - \lambda)u_i}{1}{M}}{\lfnorm{u_i}{2}{M}} \rightarrow 0$ as $i \rightarrow \infty$;
        \item $u_i \rightarrow 0$ weakly; and
        \item $\partial(\supp(u_i))$ is a $C^\infty$ $(n-1)$-submanifold of $M$.\\
    \end{enumerate}

    The construction of such a sequence will be done as in \cite{lu-zhou}.

    Given $p \in M$, let $\rtilde$ be the $C^\infty$ approximation for $r$ we have constructed in section \ref{rtil_teorema1}. Let $x, y, R$, be such that $0 < R < x < y$ and whose values will be chosen posteriorly. Let $\psi: \R \rightarrow \R$ be a cut-off function satisfying
    \[\psi(r) = \branch{1,&r\in[x/R, y/R]\\0,& r \notin [x/R-1, y/R+1]}\]
    so that $\module{\psi}$ and $\module{\psi'}$ are bounded.\\

    Given $\lambda > 0$, we define
    \[\phi(x) = \psi\left(\frac{\rtilde(x)}{R}\right)e^{i\sqrt{\lambda}\rtilde}.\]

    Hence, $\module{\phi} \leq 1$ and
    \begin{eqnarray*}
        \Delta_f\phi+\lambda\phi &=& \left[\frac1R\left(\frac{\psi''}R + 2i\sqrt{\lambda}\psi'\right)\module{\grad \rtilde}^2 +
                                \left(i\sqrt{\lambda}\psi + \frac{\psi'}{R}\right)\Delta_f \rtilde\right] e^{i\sqrt{\lambda}\rtilde}\\
                                && - \lambda\left(1 - \module{\grad \rtilde}^2\right).
    \end{eqnarray*}

    So we have
    \[\module{\Delta_f \phi + \lambda\phi} \leq \frac{C}{R} + C\module{\Delta_f\rtilde} + C\module{\grad \rtilde - \grad r}\]
    and therefore,
    \begin{eqnarray}
       \nonumber \lfnorm{\Delta_f \phi + \lambda\phi}{1}{M}
            &\leq&
        \frac{C}{R}( \widetilde{V}_f(y + R) - \widetilde{V}_f(x - R)) + C\delta(x - R)\\
        \label{desig_comum_aos_dois_casos_na_prova_do_teo1} && + C\int_{B_p(y+R)\setminus B_p(x-R)}\module{\Delta_f \rtilde}\;e^{-f}d\vol_g.
    \end{eqnarray}

    Now we divide the proof in two cases.

    First we suppose $\vol_f(M) = \infty$. By Lemma \ref{estimativa_vol_infty_teorema1}, we choose $y$ large enough such that
    \begin{equation}\label{teo1_ineq1}
        \int_{B_p(y+R)\setminus B_p(x - R)} \module{\Delta_f \rtilde}\;e^{-f}d\vol_g \leq \epsilon V_f(y + R + 1).
    \end{equation}

    On the other hand, fixing $x$ and $R$ with $x-R$ large enough, we have
    \begin{eqnarray}
        \nonumber
        \frac{C}{R}(\widetilde{V}_f(y + R) - \widetilde{V}_f(x - R))
            &<&
                \epsilon ( \widetilde{V}_f(y + R) - \widetilde{V}_f(x - R))\\
        \nonumber
            &<& \epsilon \widetilde{V}_f(y + R + 1)\\
        \label{teo1_ineq2}
            &<& 2\epsilon V_f(y + R + 1)
    \end{eqnarray}
    (in the last inequality above we have used Lemma \ref{comparacao_v_r_e_v_r_til}) and
    \begin{equation}\label{teo1_ineq3}
        C\delta(x - R) < \epsilon V_f(y + R + 1).
    \end{equation}

    By using (\ref{teo1_ineq1}), (\ref{teo1_ineq2}) and (\ref{teo1_ineq3}) in (\ref{desig_comum_aos_dois_casos_na_prova_do_teo1}), we get        \[\lfnorm{\Delta_f\phi+\lambda\phi}{1}{M} < 4\epsilon V_f(y + R + 1).\]

    Since $\phi \equiv 1$ in $B_p(y) \setminus B_p(x)$, we have $\lfnorm{\phi}{2}{M}^2 \geq \widetilde{V}_f(y) - \widetilde{V}_f(x)$. \linebreak Fixing $x$ and taking $y$ large enough, $\lfnorm{\phi}{2}{M}^2 \geq \frac{1}{2} \widetilde{V}_f(y)$, and using Lemma \ref{comparacao_v_r_e_v_r_til}, $\lfnorm{\phi}{2}{M}^2 \geq \frac{1}{4} V_f(y)$. By Lemma \ref{controle_crescimento_volume_teorema1ext}, taking $y$ even larger we have
    \linebreak $V_f(y + R + 1) \leq 2V_f(y)$, so,
    \[\displaystyle \frac{\lnorm{\phi}{\infty}{M}\;\lfnorm{(-\Delta_f - \lambda)\phi}{1}{M}}{\lfnorm{\phi}{2}{M}} < \frac{4\epsilon V_f(y + R + 1)}{\frac{1}{4} V_f(y)} \leq \frac{8\epsilon V_f(y)}{\frac{1}{4}V_f(y)} < 32\epsilon.\]

    Thus, since $\lfnorm{\phi}{\infty}{M} = 1$, we have constructed a compactly supported function $\phi$ such that
    \[\displaystyle \frac{\lnorm{\phi}{\infty}{M}\;\lfnorm{(-\Delta_f - \lambda)\phi}{1}{M}}{\lfnorm{\phi}{2}{M}} < 32\epsilon.\]
    for an arbitrarily small $\epsilon > 0$.

    In order to obtain the desired sequence $u_n$, we consider, in the above \linebreak construction, $\epsilon = 1/n$ and $u_n = \phi$, with $x$ greater than the o value of $y + R$ of the construction of $u_{n-1}$.\\

    We will now prove the case where $\vol_f(M) < \infty$.\\
    In the construction of $\rtilde$ in Proposition \ref{rtil_teorema1}, we take $\delta$ such that
    \begin{equation}\label{teo1_ext_escolha_delta}
        \delta(r) \leq \frac{1}{r}(\vol_f(M) - V_f(r)).
    \end{equation}

    By Lemma \ref{estimativa_vol_finito_teorema1ext}, we can choose $x$ large enough so that
    \begin{eqnarray*}
        \int_{B_p(x+R)\setminus B_p(x-R)} \module{\Delta_f \rtilde} \;e^{-f}d\vol_g
        &\leq& \int_{M\setminus B_p(x-R)} \module{\Delta_f \rtilde} \;e^{-f}d\vol_g\\
        &\leq& \epsilon (\vol_f(M) - V_f(x-R)) \\ &&+ 2\vol_f(\partial B_p(x - R)).
    \end{eqnarray*}

    By using this in (\ref{desig_comum_aos_dois_casos_na_prova_do_teo1}), we have
    \begin{eqnarray*}
        \lfnorm{\Delta_f\phi+\lambda\phi}{1}{M}
            &\leq&
        \frac{C}{R}( \widetilde{V}_f(y + R) - \widetilde{V}_f(x - R)) + C\delta(x - R) + \\
        &&+\;C\epsilon (\vol_f(M) - V_f(x-R)) + 2C\;\vol_f(\partial B_p(x - R))\\
            &\leq&
        \frac{C}{R}( \vol_f(M) - \widetilde{V}_f(x - R)) + C\delta(x - R) + \\
        &&+\;C\epsilon (\vol_f(M) - V_f(x-R)) + 2C\;\vol_f(\partial B_p(x - R))\\
            &\leq&
        2\frac{C}{R}( \vol_f(M) - V_f(x - R)) + C\delta(x - R) + \\
        &&+\;C\epsilon (\vol_f(M) - V_f(x-R)) + 2C\;\vol_f(\partial B_p(x - R))\\
            &\leq&
        C\left(\frac{2}{R} + \epsilon\right)(\vol_f(M) - V_f(x - R)) + C\delta(x - R) +\\
        && + 2C\;\vol_f(\partial B_p(x - R))\\
    \end{eqnarray*}
    (in the third inequality above we have used Lemma \ref{comparacao_v_r_e_v_r_til}).

    Since $M$ has finite $f$-volume, and by (\ref{teo1_ext_escolha_delta}), we may choose $R$ and $x$ so that
    \[\frac{2}{R} \leq \epsilon,\]
    and
    \[\delta(x - R) \leq \epsilon (\vol_f(M) - V_f(x - R)).\]

    Thus
    \[\lfnorm{\Delta_f\phi+\lambda\phi}{1}{M} \leq  3C\epsilon(\vol_f(M) - V_f(x - R)) + 2C\;\vol_f(\partial B_p(x - R)).\]

    By Lemma \ref{estimativa_vol_finito2_teorema1ext}, taking $x$ even larger,
    \[ 3\epsilon(\vol_f(M) - V_f(x - R)) + 2\;\vol_f(\partial B_p(x - R)) \leq 6\epsilon\;(\vol_f(M) - V_f(x)),\]
    therefore,
    \begin{equation}\label{parcial_volume_finito}
        \lfnorm{\Delta_f\phi+\lambda\phi}{1}{M} \leq 6\epsilon\;C(\vol_f(M) - V_f(x)).
    \end{equation}

    Using again the finite $f$-volume of $M$, making $y$ large enough, we get
    \[\widetilde{V}_f(y) - \widetilde{V}_f(x) \geq \frac{1}{2}(\vol_f(M) - \widetilde{V}_f(x)) \geq \frac{1}{4}(\vol_f(M) - V_f(x)),\]
    and, by (\ref{parcial_volume_finito}), we have
    \begin{eqnarray*}
        \lfnorm{\Delta_f\phi+\lambda\phi}{1}{M}
            &\leq&
        6\epsilon\;C(\vol_f(M) - V_f(x))\\
            &\leq&
        24\epsilon\;C(V_f(y) - V_f(x))\\
            &\leq&
        24\epsilon\;\lfnorm{\phi}{2}{M}^2.
    \end{eqnarray*}

    The proof follows now the case of infinite $f$-volume.

    \begin{remark}
        In \cite{charalambous2}, N. Charalambous and Z. Lu proved that the essential spectrum of $\Delta_f$ is $[0, \infty)$ under the assumptions that $\ric_f^q > -\delta^2(r)$, \linebreak $\lim_{r\rightarrow \infty} \delta(r) = 0$ and the $f$-volume of $M$ volume does not decay exponentially if it is finite. Here $\ric_f^q = \ric_f - \frac{1}{q}d f \otimes d f$ and $q > 0$.
    \end{remark}
    \begin{remark}
        In a previous version of this paper, we proved a particular case of Theorem \ref{principal1ext}, assuming the stronger hypothesis $|f| < k$. In this case, G. Wei \cite{wei1} proved that $Cr \leq \vol_f(B_p(r)) \leq Cr^{n+4k}$. These estimates make the proof of the particular case simpler, since the case $\vol_f(M) < \infty$ never happens and the conclusion of Theorem \ref{controle_crescimento_teorema1ext} follows automatically.
    \end{remark}

\section{Estimates for $\Delta r$ and $\vol (B_p(r))$ and the proof of Theorem \ref{principal2}}\label{prova_do_teorema2}

    In this Section, we suppose that $\ric_f \geq 1/2$ and $|\grad f|^2 < f$ and present some estimates for $\Delta_r$ and $V(r):=\vol(B_p(r))$. Here, we consider the approximation $\rtilde$ of $r:=d(p,\cdot)$ of the end of section \ref{rtil_teorema1}.

    \begin{lemma}\label{controle_de_lap_r_teorema2}
        If $\ric_f(\grad r,  \grad r) \geq \frac{1}{2}$ and $|\grad f|^2 \leq f$, then
        \[\lim_{r\rightarrow \infty} \Delta r \leq 0.\]
    \end{lemma}

    \proof
        Let $\gamma(t)$ be a minimizing geodesic from $p_0$. Denoting $f(t) = f(\gamma(t))$, we have
        \begin{multline*}
            0 = \frac 12 \Delta |\grad r|^2 =  |\hess r|^2 + \inner{\grad \Delta r, \grad r} + \ric(\grad r, \grad r) =\\= (\Delta r)' +  |\hess r|^2  + \ric(\grad r, \grad r),
        \end{multline*}
        so
        \begin{eqnarray*}
            (\Delta r)' &=& - |\hess r|^2 - \ric(\grad r, \grad r)\\
                      &=& - |\hess r|^2 - \ric_f(\grad r, \grad r) + \hess f(\grad r, \grad r)\\
                    &\leq& -\frac{(\Delta r)^2}{n-1} - \frac 12 |\grad r|^2 + f''\\
                      &=& -\frac{(\Delta r)^2}{n-1} - \frac 12 + f''.
        \end{eqnarray*}

        Therefore,
        \[
            \Delta r(r) - \Delta r(t_0) \leq -\frac{1}{n-1}\int_{t_0}^t (\Delta r)^2 dt - \frac{t}{2} + \frac{t_0}{2} + f'(t) - f'(t_0).
        \]

        By \cite[Proposition 4.2]{munteanu}, there exists $c$ such that \[\frac{1}{4}(t - c)^2 \leq f \leq \frac{1}{4}(t + c)^2,\]
        thus, taking the lowest positive number $t_0=t_0(\gamma)$ such that $f'(t_0) \geq 0$ and using $|f'(t)| \leq |\grad f|(\gamma(t)) \leq \sqrt{f(t)} \leq \frac{1}{2}(t + c)$, we have
        \begin{eqnarray*}
            \Delta r(t) &\leq& -\frac{1}{n-1}\int_{t_0}^t (\Delta r)^2 dt - \frac{t}{2} + \frac{t_0}{2} + \frac{1}{2}(t + c) + \Delta r(t_0) \\
                                    &=& -\frac{1}{n-1}\int_{t_0}^t (\Delta r)^2 dt + C_0,
        \end{eqnarray*}
        where \[C_0 = C_0(\gamma) = \frac{t_0}{2} + \frac{c}{2} + \Delta r(t_0).\]

        Taking $t_1$ such that $\frac{1}{4}(t_1 - c)^2 > f(p) = f(0)$, we have $f(t_1) > f(0)$. Therefore, $t_0(\gamma) \leq t_1$ for all minimizing geodesic $\gamma$ and, setting $C_1 = \sup_{\gamma}C_0(\gamma)$, we get
        \[\Delta r(t) \leq -\frac{1}{n-1}\int_{t_1}^t (\Delta r)^2 dt + C_1,\]
        which proves that $\Delta r(t)$ is a decreasing function. Given $\epsilon > 0$, taking $R > \epsilon^{-1} + C_1\epsilon^{-2} + t_1$ we have $\Delta r (x) \leq \epsilon$ for all $x \in M\setminus B_p(R)$.

    \proofend

    By \cite[Theorem 1.4]{munteanu}, we have the following volume estimates.

    \begin{lemma}
        In the hypotheses of Theorem \ref{principal2},
        \begin{enumerate}[label={(\arabic*)}]
            \item $\vol(B_p(r))$ has Euclidean volume growth, i.e., there exists $C > 0$ such that \[\vol(B_p(r)) \leq C_0 r^n.\]
            \item $\vol(B_p(r))$ has at least linear volume growth, i.e., there exists $c > 0$ such that \[\vol(B_p(r)) \geq c\,r.\]
        \end{enumerate}
    \end{lemma}

    Notice that, in particular, the conclusion (2) in the above Lemma ensures that $\vol(M) = \infty$. This, together with Lemma \ref{controle_de_lap_r_teorema2}, leads us to the following results, whose proofs are analogous to Lemmas \ref{controle_crescimento_volume_teorema1ext} and \ref{estimativa_vol_infty_teorema1}.

    \begin{lemma}\label{controle_crescimento_volume_teorema2}
        In the hypotheses of Theorem \ref{principal2}, for $R > 0$, there exists a sequence $r_k \rightarrow \infty$ of real numbers such that $\vol(B_p(r_k + R + 1)) \leq 2\vol(B_p(r_k))$.
    \end{lemma}

    \begin{lemma}\label{estimativa_vol_infty_teorema2}
        In the hypotheses of Theorem \ref{principal2}, for all $\epsilon > 0$ and $t_0$ large enough, there exists $R = R(\epsilon, t_0)$ such that, for $t > R$,
        \[\int_{B_p(t)\setminus B_p(t_0)} \module{\Delta \rtilde}d\vol_g \leq \epsilon \vol(B_p(t + 1)).\]
    \end{lemma}

    The Theorem \ref{principal2} can now be proved by constructing the same sequence of functions of the proof of Theorem \ref{principal1ext}, replacing Lemmas \ref{controle_crescimento_volume_teorema1ext} and \ref{estimativa_vol_infty_teorema1} by their corresponding versions in the new hypotheses, Lemmas \ref{controle_crescimento_volume_teorema2} and \ref{estimativa_vol_infty_teorema2} above.


\bibliography{bibliografia.txt}{}
\bibliographystyle{plain}

\end{document}